\documentclass[12pt,a4paper,notitlepage]{article}
\usepackage{amsmath,amssymb,graphicx,fullpage,authblk}

\begin{document}

\title{Recursive multivariate derivatives of $e^{f(x_1,\dots,x_n)}$ of arbitrary order}
\author{Filippo M. Miatto}
\affil{\small \emph{Institut Polytechnique de Paris}}
\affil{\small \emph{T\'el\'ecom Paris, LTCI, 19 Place Marguerite Perey, Palaiseau 91120, France}}
\date{\today}
\maketitle

\begin{abstract}
    High-order derivatives of nested functions of a single variable can be computed with the celebrated Fa\`a di Bruno's formula.
    Although generalizations of such formula to multiple variables exist, their combinatorial nature generates an explosion of factors, and when the order of the 
    derivatives is high, it becomes very challenging to compute them. A solution is to reuse what has already been computed, which is a built-in feature of recursive implementations.    
    Thanks to this, recursive formulas can play an important role in Machine Learning applications, in particular for Automatic Differentiation.
    In this manuscript we provide a recursive formula to compute multivariate derivatives of arbitrary order of $e^{f(x_1,\dots,x_n)}$ with respect to the variables $x_i$.
    We note that this method could also be beneficial in cases where the high-order derivatives of a function $f(x_1,\dots,x_n)$ are hard to compute, but where the derivatives of $\log(f(x_1,\dots,x_n))$ are simpler.
\end{abstract}

\section{introduction}
High-order derivatives of two nested functions of a single variable can be written in a compact manner by using the celebrated Fa\`a di Bruno's formula \cite{johnson2002curious}:
\begin{align}
\frac{d^n}{dx^n} g(f(x)) = n!\sum\frac{g^{(m_1+\dots+m_n)}(f(x))}{m_1!\dots m_n!}\prod_{j=1}^n\left(\frac{f^{(j)}(x)}{j!}\right)^{m_j}
\end{align}
where the sum is over all $n$-tuples of integers $m_1,\dots,m_n$ satisfying $\sum_{j=1}^n jm_j = n$. Note that we indicate the $k$-th derivative of a function $f(x)$ by $f^{(k)}(x)$.

Riordan's formula is a more convenient expression that uses the incomplete Bell polynomials (defined in the next section) \cite{johnson2002curious}:
\begin{align}
\frac{d^n}{dx^n} g(f(x)) = \sum_{k=1}^n g^{(k)}(f(x)) B_{n,k}(f^{(1)}(x),\dots,f^{(n-k+1)}(x))
\end{align}

Although generalizations of Fa\`a di Bruno's formulas for the multivariate case exist \cite{schumann2019multivariate}, they are combinatorial in nature, and don't lend themselves well to computational implementations.

A particular case of interest is when the outer function $g$ is the exponential function.
In this case, the derivatives of $g$ are all equal and can drop out of the summation.
Nevertheless, what remains is a polynomial expression with a number of terms that grows exponentially in the order of the derivatives:


\begin{align}
    \frac{\partial^{k_1+\dots+k_n}}{\partial x_1^{k_1}\dots\partial x_n^{k_n}}e^{f(x_1,\dots,x_n)} = e^{f(x_1,\dots,x_n)} Y_{k_1\dots k_n}
\end{align}

where $Y_{k_1\dots k_n}$ is a polynomial in the partial derivatives of $f$. In this work we give a recursive definition of such polynomial.

\section{Preliminaries}
\subsection{Bell Polynomials}
We begin by recalling the definition of the complete exponential Bell polynomials \cite{mihoubi2008bell}:

\begin{align}
    Y_n(z_1\dots,z_n) = \sum_{k=0}^n B_{n,k}(z_1,\dots,z_{n-k+1})
\end{align}
where $B_{n,k}(z_1,\dots,z_{n-k+1})$ are the incomplete Bell polynomials which are defined as \cite{mihoubi2008bell}

\begin{align}
    B_{n,k}(z_1,\dots,z_{n-k+1}) = \frac{1}{k!}\sum_{j_1+\dots+j_k=m\atop j_i\geq1}\frac{m!}{j_1!\dots j_m!}z_{j_1}\dots z_{j_k}
\end{align}

Already at this stage, we can simplify considerably the computation of the value of $Y_n$ at the point $(z_1,\dots,z_n)$ by using the recursive definition \cite{connon2010various}:
\begin{align}
    Y_n(z_1\dots,z_n)=\sum_{k=0}^{n-1}\binom{n-1}{k}Y_{k}(z_1,\dots,z_k)z_{n-k}
\end{align}
Together with the initial condition $Y_0 = 1$, this formula yields the complete exponential Bell polynomials of order $n$ in the variables $z_1,\dots,z_n$.

\subsection{Single variable formula}
The complete exponential Bell polynomials appear in the definition of the high-order derivatives of $e^{f(x)}$:

\begin{align}
    \frac{\partial^k}{\partial x^k}e^{f(x)} = e^{f(x)}Y_k
\end{align}
where the Bell polynomials are defined on the derivatives of $f$: $Y_k := Y_k(f^{(1)}(x),\dots,f^{(k)}(x))$.
This is a well-known result \cite{johnson2002curious} that we wish to generalize.

\section{New results}
\subsection{Recursive formula for two variables}
We first present the two variables version of our general recursive formula, as the proof is considerably simpler (the proof for $n$ variables is derived by applying the Leibniz rule as many times as necessary):
\begin{align}
    \frac{\partial^{k_1+k_2}}{\partial x_1^{k_1}\partial x_2^{k_2}}e^{f(x_1,x_2)} = e^{f(x_1,x_2)}Y_{k_1k_2}
    \label{twovariables}
\end{align}
where 
\begin{align}
    Y_{k_1k_2}=\sum_{j_1=0}^{k_1}\sum_{j_2=0}^{k_2-1}\binom{k_1}{j_1}\binom{k_2-1}{j_2}Y_{j_1j_2}f^{(k_1-j_1,k_2-j_2)}(x_1,x_2)
\end{align}

where the initial condition for the new index is $Y_{j_10}=Y_{j_1}=Y_{j_1}(f^{(1,0)}(x_1,x_2),\dots,f^{(k,0)}(x_1,x_2))$.

To reach this result, we apply the Leibnitz rule twice:

\begin{align}
    Y_{k_1k_2} &= e^{-f(x_1,x_2)}\frac{\partial^{k_1+k_2}}{\partial x_1^{k_1}\partial x_2^{k_2}}e^{f(x_1,x_2)}\\
    &=e^{-f(x_1,x_2)}\frac{\partial^{k_1}}{\partial x_1^{k_1}}\left(\frac{\partial^{k_2-1}}{\partial x_2^{k_2-1}}f^{(0,1)}(x_1,x_2)e^{f(x_1,x_2)}\right)\\
    &=e^{-f(x_1,x_2)}\frac{\partial^{k_1}}{\partial x_1^{k_1}}\sum_{s=0}^{k_2-1}\binom{k_2-1}{j_2}f^{(0,k_2-j_2)}(x_1,x_2)\frac{\partial^{j_2}}{\partial x_2^{j_2}}e^{f(x_1,x_2)}\\
    &=\sum_{j_2=0}^{k_2-1}\binom{k_2-1}{j_2}\sum_{j_1=0}^{k_1}\binom{k_1}{j_1}f^{(k_1-j_1,k_2-j_2)}(x_1,x_2)\underbrace{e^{-f(x_1,x_2)}\frac{\partial^{j_1+j_2}}{\partial x_1^{j_1}\partial x_2^{j_2}}e^{f(x_1,x_2)}}_{Y_{j_1j_2}}
\end{align}

\subsection{Recursive formula for $n$ variables (main result)}

The generalization to $n$ variables is straightforward:

\begin{align}
    \frac{\partial^{k_1+\dots+k_n}}{\partial x_1^{k_1}\dots\partial x_n^{k_n}}e^{f(x_1,\dots,x_n)} = e^{f(x_1,\dots,x_n)}Y_{k_1\dots k_n}
\end{align}
where 
\begin{align}
    Y_{k_1\dots k_n} = \sum_{j_1}^{k_1}\dots\sum_{j_n}^{k_n-1}\binom{k_1}{j_1}\dots \binom{k_n-1}{j_n}Y_{j_1\dots j_n}f^{(k_1-j_1,\dots,k_n-j_n)}(x_1,\dots,x_n)
    \label{main}
\end{align}

and the initial condition for the rightmost index is $Y_{k_1\dots k_{n-1}0}=Y_{k_1\dots k_{n-1}}$. We remark that only the last index ($j_n$ here) has an upper limit with a -1 factor.
Eq.~\eqref{main} constitutes a new result (to the best of the author's knowledge).
This is effectively a way to \emph{accumulate} derivatives of arbitrary order of the exponential of a multivariate function.

The presence of the exponential is not strictly necessary, as our recursive method can also be restated as:
\begin{align}
    \frac{\partial^{k_1+\dots+k_n}}{\partial x_1^{k_1}\dots\partial x_n^{k_n}}g(x_1,\dots,x_n) = g(x_1,\dots,x_n)T_{k_1\dots k_n}
\end{align}
where 
\begin{align}
    T_{k_1\dots k_n} = \sum_{j_1}^{k_1}\dots\sum_{j_n}^{k_n-1}\binom{k_1}{j_1}\dots \binom{k_n-1}{j_n}T_{j_1\dots j_n}\frac{\partial^{j_1+\dots+j_n}}{\partial x_1^{j_1}\dots\partial x_n^{j_n}}\log(g(x_1,\dots,x_n))
\end{align}

\section{Complexity and Benchmarks}
The time complexity of the whole procedure is dominated by the complexity of computing the derivatives of $f$.
If we don't exploit any structure of $f$ or reuse computed values wherever possible, the time complexity of computing the derivative of order $(k_1\dots k_n)$ of $e^{f(x_1,\dots,x_n)}$ requires evaluating the derivatives of $f$ $\prod_{j=1}^n k_j$ times.

It's important to notice that if the derivatives of $f$ are evaluated numerically, or if they are evaluated at a point, $Y_{j_1\dots j_n}$ at every step is just a number.
In particular, if we have a closed form for the derivatives of $f$, this process can be extremely efficient, as it just requires evaluating a \emph{function} (which can be compiled) $\prod_{j=1}^n k_j$ times.

We now compare our method, written in Mathematica, against the standard Mathematica implementation for two select example functions.
To assess the efficiency of both implementations fairly, we clear the memoization cache after each value is computed.
This means that both implementations begin from scratch for each point of the graph. Note that in neither of the implementations we take into account the particular nature of $f$.

For the first problem, the function $f$ at the exponent is polynomial of four variables:
\begin{align} 
    F(x_1,x_2,x_3,x_4) = \exp(x_1x_2x_3x_4 + x_1^2x_2^2x_3^2x_4^2 + x_1^3x_2^3x_3^3x_4^3)
\end{align}

and we evaluate the time to compute the following high-order derivatives for $0\leq k \leq 6$:
\begin{align} 
    \frac{\partial^{12+k}}{\partial^kx_1\partial^4x_2\partial^4x_3\partial^4x_4}F(x_1,x_2,x_3,x_4)\Bigr|_{x_i=1}
\end{align}

\begin{figure}[ht]
    \begin{center}
        \includegraphics[scale=0.8]{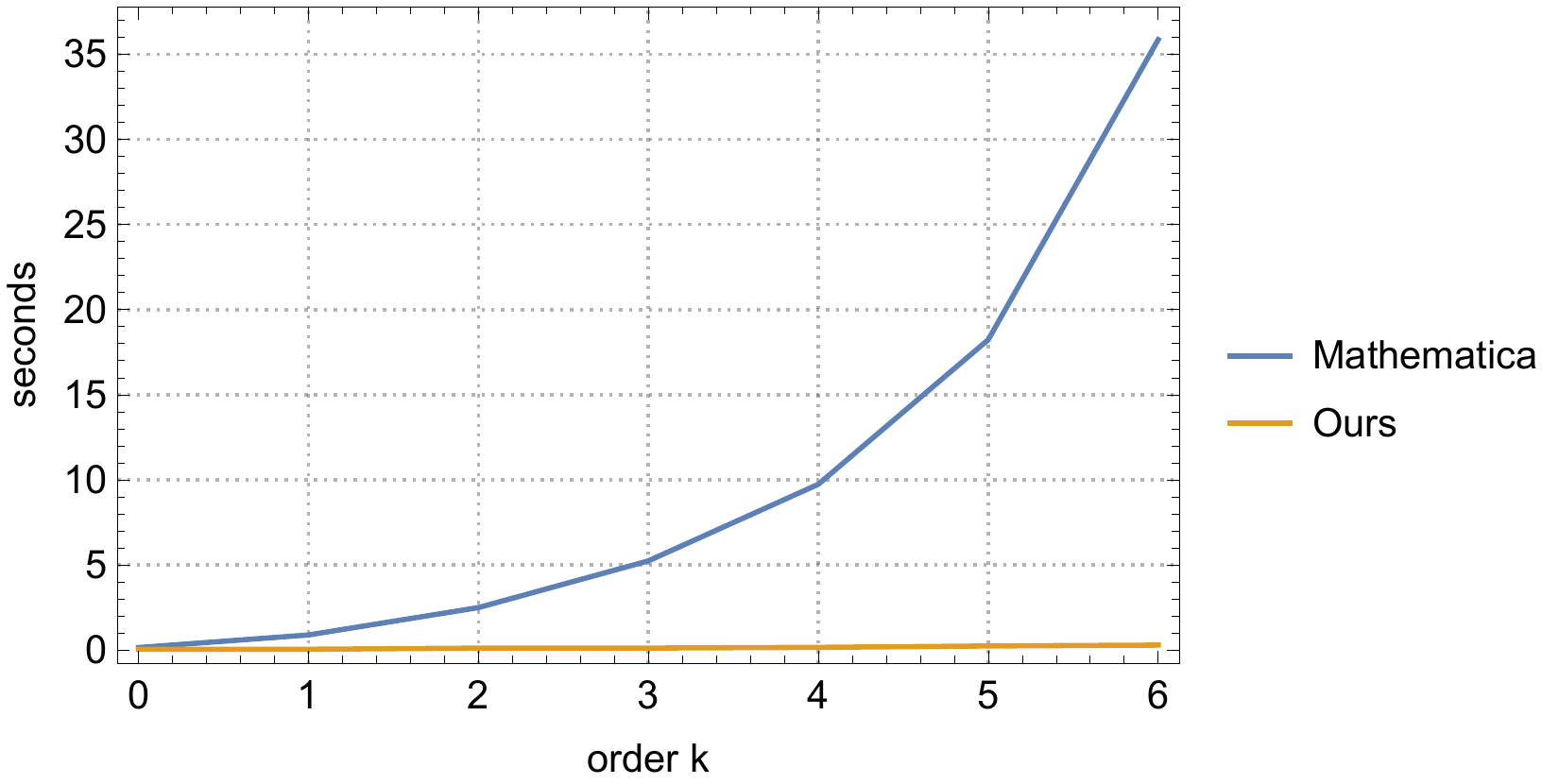} 
        \caption{Time to compute high-order derivatives of the exponential of a polynomial function. The longest time for our implementation (for $k=6$) is about 0.3s.}
    \end{center}
\end{figure}

For the second problem, we choose a non-polynomial function:
\begin{align} 
    G(x_1,x_2,x_3,x_4) = \exp(x_1x_2x_3\sin(x_4) + x_1x_2\sin(x_3)x_4 + x_1\sin(x_2)x_3x_4 + \sin(x_1)x_2x_3x_4)
\end{align}

and again we evaluate the time to compute the following high-order derivatives for $0\leq k \leq 6$:
\begin{align} 
    \frac{\partial^{12+k}}{\partial^kx_1\partial^4x_2\partial^4x_3\partial^4x_4}G(x_1,x_2,x_3,x_4)\Bigr|_{x_i=1}
\end{align}

\begin{figure}[ht!]
    \begin{center}
        \includegraphics[scale=0.8]{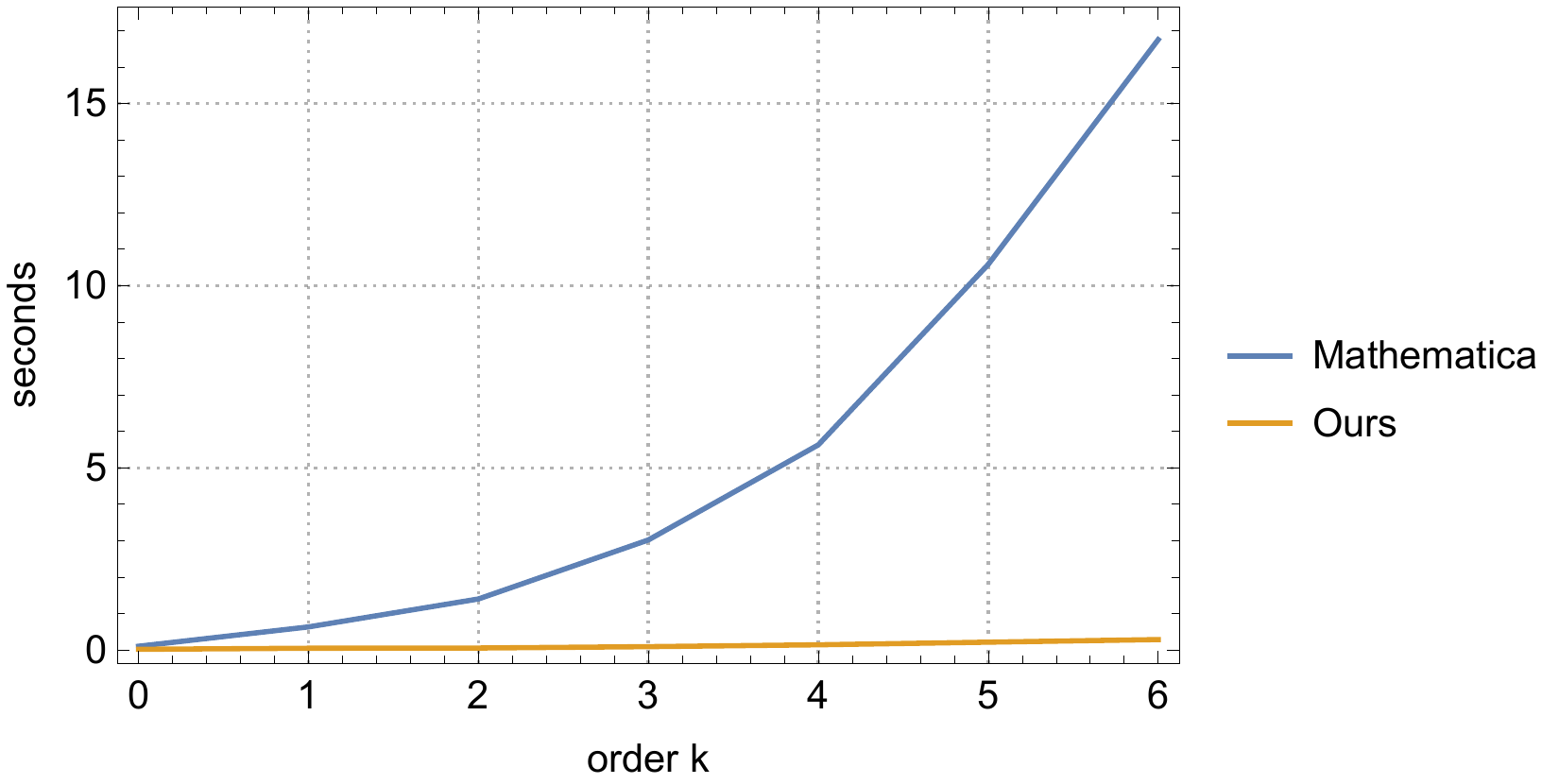} 
        \caption{Time to compute high-order derivatives of the exponential of a non-polynomial function. The longest time for our implementation (for $k=6$) is about 0.2s.}
    \end{center}
\end{figure}

The difference in performance is more striking the more variables are involved.

As pointed out above, the standard Mathematica implementation could in principle be rewritten to take into account the specific properties of the functions being exponentiated.
However, this would require a special implementation for each case, whereas we used a general implementation using the standard language features.
We shoud also note that both implementations return matching values, up to machine precision.

\section{Conclusions}
We have presented an efficient recursive algorithm for computing derivatives of arbitrary order of the exponential of a multivariate function.
Our implementation can outperform significantly other state of the art solutions, such as what one could achieve with the standard methods in Mathematica.

Our hope is that this new method will be implemented in symbolic and automatic differentiation software, as well as being used in formal studies in calculus, functional analysis and related fields.

\section{Acknowledgements}
The author thanks Markus Scheuer for help with the proof of \eqref{twovariables}.

\bibliographystyle{unsrt}
\bibliography{expbiblio}

\end{document}